\documentstyle[12pt]{article}
\title{A Noether-Lefschetz Theorem and applications}
\author{ Kirti Joshi\\
	School of Mathematics,\\
	Tata Institute of Fundamental Research,\\
	Homi Bhabha Road,
	Bombay 400 005.\\
	E-mail: Orion@tifrvax.bitnet
	}
\date{Preliminary version}
\begin{document}
\maketitle

\newtheorem{thm}{Theorem}[section]
\newtheorem{propose}{Proposition}[section]
\newtheorem{lemma}[propose]{Lemma}
\newtheorem{cor}[propose]{Corollary}
\newtheorem{sublem}[propose]{Sublemma}
\newtheorem{schol}{Scholium}
\newtheorem{problem}{Problem}
\newtheorem{conj}{Conjecture}
\newtheorem{question}{Question}
\newtheorem{thmp}{Theorem}
\renewcommand{\thethmp}{{\thethm${}'$}}
\newtheorem{notation}{Notation}

\renewcommand{\P}{{\ifmmode{\bf P}\else{\bf P}\fi}}   
\newcommand{\Q}{\mbox{$\bf Q$}}     
\newcommand{\A}{\mbox{$\bf A$}}     
\newcommand{\R}{\mbox{$\bf R$}}     
\newcommand{\C}{{\ifmmode{\bf C}\else{\bf C}\fi}}     
\newcommand{\Z}{\mbox{$\bf Z$}}     
\newcommand{\N}{\mbox{$\bf N$}}     
\newcommand{\G}{\mbox{$\bf G$}}     
\newcommand{\implies}{\mbox{$\Rightarrow$}}
\newcommand{\ie}{{\it i.e.\/},\ }
\newenvironment{proof}{
                       \trivlist \item[\hskip \labelsep{\bf Proof}:]
                      }{
			\hfill$\Box$\endtrivlist
		      }

\newenvironment{defn}{
			\trivlist \item[\hskip \labelsep{\bf Definition}:]
		      }{
			\endtrivlist}
\newenvironment{claim}{
			\trivlist \item[\hskip \labelsep{\bf Claim}:]
		      }{
			\endtrivlist}
\newenvironment{rmk}{
			\trivlist \item[\hskip \labelsep{\bf Remark}:]
		      }{
		     	\endtrivlist}

\let\union=\cup
\let\intn=\cap
\let\tensor=\otimes
\let\into=\hookrightarrow
\let\tilde=\widetilde
\let\hat=\widehat
\let\longto=\longrightarrow
\let\lra=\longrightarrow
\newcommand{\Spec}{\mathop{\rm Spec}\nolimits}    
\newcommand{\Pic}{\mathop{\rm Pic}\nolimits}      
\newcommand{\codim}{\mathop{\rm codim}\nolimits}  
\newcommand{\coker}{\mathop{\rm coker}\nolimits}  %
\newcommand{\im}{\mathop{\rm im}\nolimits}        
\let\image=\im
\newcommand{\Tor}{\mathop{\rm Tor}\nolimits}
\def\E{{\cal E}}
\def\O{{\cal O}}
\let\da=\downarrow
\def\vnl{\mathop{\rm VNL}\nolimits}
\def\inl{\mathop{\rm INL}\nolimits}
\def\rk{{\mathop{\rm rk}\nolimits}}
\def\oxy{\Omega_{(\X,\Y)}}
\def\tX{{\tilde X}}
\def\tY{{\tilde Y}}
\def\X{{\cal X}}
\def\Y{{\cal Y}}
\def\I{{\cal I}}
\def\bQ{\bar{\bf Q}}
\let\into=\hookrightarrow
\let\isom=\simeq
\let\la=\longrightarrow
\let\abuts=\implies
\def\bull{{\scriptscriptstyle\bullet}}
\def\reg{\mathop{\rm reg}\nolimits}
\newcommand{\etab}{{\overline\eta}}
\def\specialref#1#2{{\rm(\ref{#1}.{#2})}}
\let\ra=\longrightarrow
\let\RA=\Longrightarrow
\def\mapright#1{\smash
                {\mathop{\longrightarrow}\limits^{#1}}
                } 
\def\mapdown#1{\Big\downarrow
                \rlap{
                        $\vcenter{\hbox{$\scriptstyle{#1}$}}$
                       }
                }  
\def\mapup#1{\Big\uparrow
		\rlap{
			$\vcenter{\hbox{$\scriptstyle{#1}$}}$
			}
		}
\def\Matrix#1{
                \def\normalbaselines{
                                \baselineskip20pt
                                \lineskip3pt\lineskiplimit3pt
                                }
                \matrix{#1}
                }

\section{Introduction}

        In this paper we generalize the classical Noether-Lefschetz
Theorem (see \cite{Lefschetz}, \cite{GH}) to arbitrary smooth
projective threefolds.  More specifically, we prove that given any
smooth projective threefold $X$ over complex numbers and a very ample
line bundle $L$ on $X$, there is an integer $n_0(X,L)$ such that if $n
\geq n_0(X,L)$ then  the Noether-Lefschetz locus of the linear
system $H^0(X,L^n)$ is a countable union of proper closed subvarieties
of $\P(H^0(X,L^n)^*)$ of codimension at least two.  In particular, the
{\em general singular member} of the linear system $H^0(X,L^n)$ is not
contained in the Noether-Lefschetz locus.  This generalizes the
results of \cite{GH}.  In \cite{Madhav}, we find a
conjecture due to M.~V.~Nori, which generalizes the Noether-Lefschetz
theorem for codimension one cycles on smooth projective threefolds to
higher codimension cycles on arbitrary smooth projective varieties.
As an application of our main theorem we prove a result which can be
thought of as a weaker version of Nori's conjecture for {\em
codimension two cycles} on smooth projective threefolds.

        The idea of the proof is borrowed from an elegant paper of
M.~Green (\cite{MarkGreen}) and the work of M.~V.~Nori (see
\cite{Madhav}). In \cite{MarkGreen} it was shown that the
classical Noether-Lefschetz theorem can be reduced to a coherent
cohomology vanishing result and the required vanishing was also proved
for $\P^3$.  Though we have used this idea, we prove the required
vanishing by combining the techniques of \cite{Madhav} and of
\cite{MarkGreen}.

        The paper is organized as follows. In the next section
(Section \ref{Generalities and notations}) we set up the basic
notations and terminology. Here we have also collected a few facts
which will be used throughout this paper. In Section~\ref{The
Noether-Lefschetz machine} we set up the basic technical machinery. In
this section, we reduce the Infinitesimal Noether-Lefschetz theorem to
certain coherent cohomology vanishing statements. This is the most
crucial part of the paper. Following the unpublished work of N. Mohan
Kumar and V. Srinivas (see \cite{Mohan-Srinivas}), in
Section~\ref{The Noether-Lefschetz Locus}, we show how the
Infinitesimal
Noether-Lefschetz theorem can be used to prove the global
Noether-Lefschetz theorem.  The global Noether-Lefschetz theorem is
thus reduced to Infinitesimal Noether-Lefschetz theorems, which in
turn are deduced from coherent cohomology vanishing results. In
Section~\ref{A General Noether-Lefschetz Theorem}, we prove our main
result, Theorem~\ref{main theorem}. By the results of the
Sections~\ref{The Noether-Lefschetz machine} and \ref{The
Noether-Lefschetz Locus}, to prove our main theorem, we are reduced to
proving several cohomology vanishing results. This program is carried
out in Section~\ref{A General Noether-Lefschetz Theorem}. The
technique of the proofs of this section are based on Green's work
(see \cite{MarkGreen}) and a modification of this method due to
Paranjape (see \cite{KP}).  Finally, in the last section we give an
application of our main theorem to codimension two cycles on smooth
projective threefolds.

        This paper could not have been written without all the help
that we have received from N. Mohan Kumar and Kapil Paranjape. They
have generously shared their ideas and insights on the problems
considered in this paper; moreover N. Mohan Kumar suggested the
problem and also explained to us his unpublished work (with V.
Srinivas) and the work of Green.  Kapil Paranjape patiently
explained to us the work of Nori and the subsequent
simplifications of Nori's work due to him.  We would like to thank
both N.~Mohan Kumar and Kapil Paranjape for all their help and
encouragement (and also their patience) without which this paper might
not have been written.

        We would also like to thank Madhav Nori and V.~Srinivas for
numerous conversations and suggestions.

\section{Generalities and notations}\label{Generalities and
notations}
        Let $X/\C$ be a smooth projective threefold, $L$ an ample line
bundle on $X$, spanned by its global sections. For any linear system $
W \subset H^0(X,L) $ we will write $S=\P(W^*)$, and for any point
$s\in S$ we will  write $T = \Spec(\O_{\P(W^*)}/m^2_s)$, where $m_s$
is the maximal ideal of the point $s$.  The base point free linear
system $H^0(X,L)$ gives a universal family of hypersurfaces with the
parameter space $\P(H^0(X,L)^*)$.  Moreover, if $W\subset H^0(X,L)$ is
any sub linear system, then we can restrict the above universal family,
to the linear subspace $\P(W^*) \subset \P(H^0(X,L)^*)$, and get a
universal family of hypersurfaces $\Y$ parametrised by $\P(W^*)$.
Further, we will also write $\X = X \times \P(W^*)$.  Note that we
have suppressed the dependence of $\X,\Y$ and $T$ on $W$ and $s$. In
our discussions, it will be clear which $W,s$ the notation refers to,
and there will be no danger of confusion.

When the linear system $W$ is base point free, we can also describe
the universal family as follows.  We denote by $M(L,W)$ the locally
free sheaf defined by the exact sequence
\begin{equation}
 0 \to M(L,W) \to W \tensor \O_X \to L \to 0.
 \end{equation}
 Further as a notational convenience, if $W=H^0(X,L)$, we will write
$M(L) = M(L,W)$. Then  $\Y = \P_X(M(L,W)^*)$. The dual of the exact
sequence above then gives an embedding $\Y \into \X$.

For any $S$-scheme $S'$, denote by $\X_{S'} = \X \times_{S} S'$. In
particular we have the ``universal infinitesimal deformations''
corresponding to a point $s\in S$, given by $\X_{T}, \Y_T$. Note that
for any $S' \to S$, we have an embedding $\Y_{S'} \into \X_{S'}$.

We have the projections $p_X : \X \to X$ and $p_S:\X \to S$. The
restrictions of these morphisms to $\Y$ will also be denoted by the
same symbols. Note that $p_S:\Y \to S$ is smooth over $s\in S$ if the
zero scheme of the section $s\in W$ is a smooth subscheme of $X$. We
will denote by $Y_s$ the zero scheme of $s$.

        We also have on $\X$ (resp. on $\Y$) the de Rham complex
on $\X$ (resp. on $\Y$) denoted by
$\Omega_\X^\bull$ (resp.  $\Omega_\Y^\bull$).
We define  a complex
on $\X$, denoted by $\Omega_{(\X,\Y)}^\bull$ by the exact
sequence
\begin{equation}
         0 \to \Omega_{(\X,\Y)}^\bull \to \Omega_{\X}^\bull \to
                \Omega_\Y^\bull  \to 0.
 \end{equation}
The first two are complexes of locally free sheaves on $\X$ and the
latter is a complex of locally free sheaves on $\Y$.
 Further as the natural projections  $\X \to X$ and $\Y \to X$ are
smooth morphisms, we see that we also have a commutative diagram:
 \label{main diagram}
\begin{equation}
        \Matrix{
  &   &    0     &   &         0        &   &      0       &    &    \cr
  &   & \da      &   & \da              &   & \da          &    &    \cr
0 &\to& p_X^*(\Omega_X^1) \tensor \O_\X(-\Y) &\to& \Omega^1_{(\X,\Y)} &\to&
 \Omega^1_{(\X,\Y)/X}&\to&0\cr
  &   & \da      &   & \da              &   & \da          &    &    \cr
0 &\to& p_X^*(\Omega^1_X) &\to& \Omega^1_{\X}      &\to& \Omega^1_{\X/X}& \to
 &
    0\cr
  &   & \da      &   & \da              &   & \da          &    &    \cr
0 &\to& p_X^*(\Omega^1_X)\big|_\Y &\to& \Omega^1_{\Y}      &\to&
\Omega^1_{\Y/X}
   & \to  & 0\cr
  &   & \da      &   & \da              &   & \da          &    &    \cr
  &   &    0
&   &         0        &   &      0       &    &    \cr
        }
 \end{equation}
        Note that as $\Y$ is a divisor in $\X$, the sheaves
$\Omega^1_{(\X,\Y)}$ and $\Omega^1_{(\X,\Y)/X}$ are vector bundles on $\X$.


\begin{defn}\label{VNL}
    We will say the Vanishing Noether-Lefschetz
condition is valid for $(X,L,W,s)$ where $W$ is a
linear system contained in $H^0(X,L)$,
and $s\in W$ if the following
assertion is valid:
            $$H^2( X , \Omega^1_{(\X_T,\Y_T)}\big|_{X\times s} ) = 0.$$
        We will symbolically denote this hypothesis by $\vnl(X,L,W,s)$.
 \end{defn}

\begin{defn}
        If $X,L$ are as above and $W\subset H^0(X,L)$ is any linear
system we will say that the Infinitesimal Noether-Lefschetz Theorem is
valid at $s\in W$ if in the commutative diagram
        $$\Matrix{
                \Pic(X) & \to & \Pic(\Y_T) \cr
                   \searrow & & \swarrow \cr
                        & \Pic(Y_s) &           \cr
                }
        $$
 the equality  $\image(\Pic(\Y_T) \to \Pic(Y_s)) = \image(\Pic(X)\to
 \Pic(Y_s) ) $ is valid.
        We will denote this hypothesis by $\inl(X,L,W,s)$.
 \end{defn}

The key observation which relates these two definitions is the
following proposition, which is just a reformulation of Remark 3.10
of \cite{Madhav}.

\begin{propose}\label{vnl gives inl}
    Let $X$ be a smooth projective threefold over complex numbers,
$L$ a very ample line bundle on $X$, $W\subset H^0(X,L)$ a base point
free linear system. If $s\in W$ cuts out a smooth divisor $Y_s$ on
$X$ and if $\vnl(X,L,W,s)$ is valid then so is $\inl(X,L,W,s)$.
 \end{propose}

\begin{proof}
    The morphism $p_S: \X \to S$ gives rise to an exact sequence
        $$ 0 \la p_S^*\Omega^1_S \la \Omega^1_\X \la \Omega^1_{\X/S}
                \la 0,$$
 and a similar exact sequence for the morphism $p_S:\Y \to S$.  Since
$Y_s$ is a smooth divisor on $X$, the morphism $p_S:\Y\to S$ is smooth
over $s\in S$. Hence we have the following commutative diagram:
 \begin{equation}
        \Matrix{
  &   &    0     &   &         0        &   &      0       &    &    \cr
  &   & \da      &   & \da              &   & \da          &    &    \cr
0 &\to& p_S^*(\Omega_S^1) \tensor \O_X(-Y_s) &\to&
                \Omega^1_{(\X,\Y)}\big|_{X_s} &\to&
                        \Omega^1_{(\X,\Y)/S}\big|_{X_s}&\to&0\cr
  &   & \da      &   & \da              &   & \da          &    &    \cr
0 &\to& p_S^*(\Omega^1_S)\big|_{X_s} &\to&
        \Omega^1_{\X}\big|_{X_s}      &\to&
                \Omega^1_{\X/S}\big|_{X_s}& \to  & 0\cr
  &   & \da      &   & \da              &   & \da          &    &    \cr
0 &\to& p_S^*(\Omega^1_S)\big|_{Y_s} &\to&
        \Omega^1_{\Y}\big|_{Y_s}      &\to&
                \Omega^1_{\Y/S}\big|_{Y_s}& \to  & 0\cr
  &   & \da      &   & \da              &   & \da          &    &    \cr
  &   &    0     &   &         0        &   &      0       &    &    \cr
        }
 \end{equation}
        This diagram gives rise to the following cohomology diagram:
        \begin{equation}\let\scty=\scriptstyle
        \Matrix{
  &      & \scty{H^1( X , \Omega_X^1)}                                & \scty{
\mapright{}}
                & \scty{H^2( X ,  p_S^*(\Omega_S^1) \big|_{X_s})}     \cr
  &      &  \scty{\mapdown{\beta}}                             &
                & \scty{\mapdown{}}      \cr
 \scty{H^1(Y_s,\Omega^1_{\Y}\big|_{Y_s})} & \scty{\mapright{}}
      & \scty{H^1(Y_s,\Omega^1_{Y_s})}                             & \scty{
\mapright{\alpha}}
                & \scty{H^2(Y_s,p_S^*(\Omega_S^1)\big|_{Y_s})} \cr
\scty{\mapdown{}}  &      & \scty{\mapdown{}}      &
                & \scty{\mapdown{}}                            \cr
        \scty{H^2(X,\Omega^1_{(\X,\Y)}\big|_{X_s})}  & \scty{\mapright{}}
                & \scty{H^2 ( X, \Omega^1_{(\X,\Y)/S}\big|_{X_s})}
                & \scty{\mapright{}}
                & \scty{H^3( X , p_S^*(\Omega^1_S)\tensor \O_X(-Y_s))}      \cr
                }
        \end{equation}
        This  diagram is obtained from the previous diagram by taking
the long exact cohomology sequence and noting that we have the
identifications:

        \begin{eqnarray*}
           \Omega^1_{\X/S}\big|_{X_s} & \isom & \Omega^1_{X_s} \\
           \Omega^1_{\Y/S}\big|_{Y_s} & \isom & \Omega^1_{Y_s}
         \end{eqnarray*}
        By hypothesis, the first term on the bottom row (from the
left) is zero (this vanishing is just $\vnl(X,L,W,s)$). Then it
follows from the injectivity of the bottom row map that $\ker(\alpha)
\subset \image(\beta)$. Further note that one has the trivial
identifications:
        \begin{eqnarray*}
           H^2(X_s,p_S^*(\Omega^1_S)\big|_{X_s}) & \isom &
                                H^2(X,\O_X)\tensor \Omega^1_S\tensor k(s) \\
           H^2(Y_s,p_S^*(\Omega^1_S)\big|_{Y_s}) & \isom &
                                H^2(Y_s,\O_{Y_s})\tensor \Omega^1_S \tensor
k(s)
   ,
         \end{eqnarray*}
 where $k(s)$ denotes the residue field of $s$. Using these
identifications, one notes that the top two rows of this diagram can
be identified with the ``Kodaira-Spencer maps'' for the infinitesimal
deformations $\X_T$ and $\Y_T$ respectively (see \cite{CGGH})
(strictly speaking this identification above should be carried
out on $T$).

        Since $\ker(\alpha) \subset \image(\beta)$, we see that any
any class $c\in H^1(Y_s,\Omega^1_{Y_s})$ which deforms infinitesimally
is the image of a class in $H^1(X,\Omega^1_X)$. Further as $Y_s$ is
the zero scheme of a section of a very ample line bundle, one has a
surjection $\Pic^0(X) \to \Pic^0(Y_s)$ (this follows easily from
Kodaira vanishing theorem). Moreover, it is easy to see, using the
exponential sequence, that the homomorphism $\Pic(X) \to \Pic(Y_s)$
has torsion free cokernel.  So that we  see at once that in the
commutative diagram
        $$
        \Matrix{
                \Pic(X) &      \to    & \Pic(\Y_T) \cr
               \searrow &             & \swarrow \cr
                        &  \Pic(Y_s)  & \cr
                        }
        $$
 any line bundle which is in the image of the map $\Pic(\Y_T) \to
\Pic(Y_s)$, comes from $\Pic(X)$. Thus we have proved the proposition.
 \end{proof}

This proposition, though a technical assertion, is a crucial result
from the point of view of this paper. In the next section we use this
result to systematically reduce Infinitesimal Noether-Lefschetz
Theorem to a coherent cohomology vanishing which can be checked in practice.

\section{The Noether-Lefschetz machine}\label{The Noether-Lefschetz machine}
        The condition $\vnl(X,L,W,s)$ is now  cast into a more
manageable form. The technique which is employed in the proof of the
following proposition is going to be applied repeatedly in different
contexts throughout this section and hence in the subsequent proofs of
this section, we will give the important points and leave the details
to the reader.

\begin{propose}\label{base free and smooth}
         Let $(X,L,W,s)$ be such that $W\subset H^0(X,L)$ is a base
point free linear system and $s\in W$ cuts out a smooth divisor on
$X$. Assume that 
the following assertions are valid:
 \begin{description}
            \item[{\rm(\thepropose.1)}] $H^1( X , \Omega^2_X \tensor L ) = 0$,
a
   nd
            \item[{\rm(\thepropose.2)}] $H^1( X ,  M(L,W) \tensor K_X \tensor L
   ) = 0 $.
\end{description}
        Then $\vnl(X,L,W,s)$ is valid.
 \end{propose}

\begin{proof}
        We closely follow the techniques of proof in \cite{Madhav}.
Indeed, we make repeated application of this method throughout this
section.

         Since the morphisms $p_X:\X \to X$ and $p_X\big|_\Y: \Y \to
X$ are smooth we have the following exact sequence of vector bundles
on $\X$:
\begin{equation}
  0 \to p_X^*(\Omega^1_X)\tensor \O_\X(-\Y) \to
                \Omega^1_{(\X,\Y)} \to \Omega^1_{(\X,\Y)/X} \to 0.
 \end{equation}
         This exact sequence is just the top row of the commutative
diagram (\ref{main diagram}). Further, note that we can restrict this
exact sequence to $X_s$ and get the exact sequence:
 \begin{equation}
         0 \to p_X^*(\Omega^1_X) \tensor \O_\X(-\Y) \big|_{X_s}
                 \to \Omega^1_{(\X,\Y)} \big|_{X_s} \to
                        \Omega^1_{(\X,\Y)/X}\big|_{X_s} \to 0.
 \end{equation}
        Thus to prove that the middle term of the above exact sequence
has no $H^2$, it suffices to prove that the extreme terms have no
$H^2$. Now as $\O_\X(-\Y)\big|_{X_s} = L^{-1}$,we see that
 \begin{equation}
 H^2(X , p_X^*(\Omega^1_X)\tensor\O_\X(-\Y)\big|_{X_s}) = H^2(X ,
                \Omega^1_X\tensor L^{-1}).
 \end{equation}
         And so the vanishing of this is just the hypothesis
\specialref{base free and smooth}{1}, by Serre duality.
Now we prove that the
second term also has no $H^2$. This is done by the following:
\begin{claim}\label{from Kapils paper}
    $$\Omega^1_{(\X,\Y)/X}\big|_{X_s} = M(L,W)^* \tensor L^{-1}.$$
 \end{claim}
        \begin{proof}
                This fact is easily verified, it also follows from an
explicit resolution of the sheaves $\Omega^i_{(\X,\Y)/X}$ constructed
in \cite{KP}.
 \end{proof}
        Now the second vanishing is just our hypothesis
\specialref{base free and smooth}{2}, after applying Serre duality.
This proves Proposition \ref{base free and smooth}.
 \end{proof}

        Before we proceed further we need some notations. If $x\in X$
is any point, we denote by $\tX_x$ the blowup of $X$ along $x$;
$\pi_x:\tX_x \to X$ the blowup morphism and $E_x$ the exceptional
divisor. We write $L'_x = \pi_x^*(L)\tensor\O_{\tX_x}(-E_x)$.
Frequently, when there is no chance of confusion, we will
suppress the subscript $x$ from the above notations.

Suppose $W\subset H^0(X, m_x \tensor L)=H^0(\tX_x,L'_x)$ is a subspace.
We assume that the linear system $W$ is base point free
on $\tX$. By Proposition \ref{base free and smooth},
$\inl(\tX,L'_x,W,s)$ is reduced to  a vanishing on $\tX$.  The next
proposition reduces the vanishing on $\tX$ to a coherent cohomology
vanishing on $X$.

\begin{propose}\label{base point and smooth}
         Let $(X,L,W,s)$ be as above. Suppose
 the following assertions are valid:
        \begin{description}
            \item[{\rm(\thepropose.1)}] $H^1( X ,
                        (\pi_x)_*(\Omega^2_{\tX_x}(-2E_x))\tensor L) =0$, and
             \item[{\rm(\thepropose.2)}] \label{b-p-s-2} $H^1( X , (\pi_x)_*
                        (M(L'_x,W))\tensor K_X \tensor  L) = 0$.
         \end{description}
        Then $\vnl(\tX_x , L'_x , W, \tilde{s})$ is valid.
 \end{propose}

\begin{proof}
        The proof is similar to the one given earlier, though it is a
bit more involved. For the purpose of the proof let us use the
following notations: we write $\X = \tX_x \times S$, $S=\P(W^*)$,
$\Y=\P_{\tX_x}(M(L'_x,W)^*)$ and $D=E_x
\times S$. Note that $D$ is a divisor on $\X$ and that $\O_\X(D) =
p_{\tX_x}^*(\O_{\tX_x}(E_x))$.

        We have the following exact sequence on $\X$.
  $$ 0 \to \Omega^1_{(\X,\Y)} \to \Omega^1_{(\X,\Y)}(D) \to
                \Omega^1_{(\X,\Y)} (D) \big|_D \to 0.$$
        Now we can restrict this exact sequence to the fibre
$\X_s=\tX \times \{\tilde s\}$ of the projection morphism $\X \to
S$ over $\tilde s\in S$. Thus we have an exact sequence
          $$ 0 \to \Omega^1_{(\X,\Y)}\big|_{\X_s} \to
                \Omega^1_{(\X,\Y)}(D)\big|_{\X_s} \to
                        \Omega^1_{(\X,\Y)} (D) \big|_D\big|_{\X_s} \to 0.$$
Thus to prove that the first term on the left has no $H^2$
(that is to say $\vnl(\tX_x , L'_x , W, \tilde{s})$), it
suffices to prove the following assertions:
\begin{description}
    \item[{\rm(\thepropose.3)}] $H^2(\X_s,
                \Omega^1_{(\X,\Y)}(D)\big|_{\X_s} ) = 0$, and
    \item[{\rm(\thepropose.4)}] $H^1( \X_s,
                \Omega^1_{(\X,\Y)}(D)\big|_{D}\big|_{\X_s} )= 0$.
 \end{description}
        These are proved in  Lemma \ref{proof of 3.2.3} and Lemma
\ref{proof of 3.2.4}.
 \end{proof}

\begin{lemma}\label{proof of 3.2.3}
    The hypotheses of Proposition \ref{base point and smooth} imply
\specialref{base point and smooth}{3}.
 \end{lemma}

\begin{proof}
    First note that $D\big|_{\X_s} = E_x$ (recall that $D=E_x\times
S$, so this equality is obvious). Now we have the exact sequence
        $$ 0 \to p_{\tX}^*(\Omega^1_{\tX}) \tensor \O_X(-\Y) \to
                \Omega^1_{(\X, \Y)} \to \Omega^1_{(\X,\Y)/\tX} \to 0.$$
Twisting this exact sequence by $\O_\X(D)$ we get
        $$\scriptstyle{ 0 \to p_{\tX_x}^*(\Omega^1_{\tX_x}) \tensor \O_X(-\Y)
                \tensor \O_\X(D) \to
                        \Omega^1_{(\X, \Y)} \tensor \O_\X(D) \to
                        \Omega^1_{(\X,\Y)/\tX_x} \tensor \O_\X(D) \to 0.}$$
Now restricting to $\X_s$, we have
 $$ 0  \to p_{\tX}^*(\Omega^1_{\tX}(E_x)) \tensor \O_\X(-\Y)\big|_{\X_s} \to
        \Omega^1_{(\X,\Y)}(D)\big|_{\X_s} \to
                \Omega^1_{(\X,\Y)/\tX}(D) \big|_{\X_s} \to 0,$$
 where we have used the fact that $\O_\X(D) =
p_{\tX}^*(\O_{\tX}(E_x))$ in the term on the left. So the middle term
has no $H^2$ if the extreme terms have no $H^2$. So we have to check
that the hypotheses of the Proposition \ref{base point and smooth}
ensure this.

        Note that $\O_\X(-\Y) \big|_{\X_s} = L'^{-1}_x$, and as
$p^*_\tX(\Omega^1_\tX (E_x))\big|_{\X_s} = \Omega^1_{\tX}(E_x)$, we
have to prove that the hypotheses of the proposition imply that
\begin{description}
        \item[{\rm(\thepropose.5)}] $H^2(\tX, \Omega^1_\tX(E_x) \tensor
                L'^{-1}_x) = 0$,  and
        \item[{\rm(\thepropose.6)}] $H^2(\tX
                ,\Omega^1_{(\X,\Y)/\tX}(D)\big|_{\X_s}) = 0 $.
 \end{description}
 We will use the Leray spectral sequence to prove that the hypothesis
\specialref{base point and smooth}{1} will imply
\specialref{proof of 3.2.3}{5} and hypothesis
\specialref{base point and smooth}{2} will imply
\specialref{proof of 3.2.3}{6}.  This will
prove Lemma \ref{proof of 3.2.3}.
 \end{proof}

\begin{sublem}\label{sublemma1}
        Hypothesis \specialref{base point and smooth}{1} implies
\specialref{proof of 3.2.3}{5}.
 \end{sublem}
\begin{proof}
        By Serre duality we are reduced to proving:
        $$H^1( \tX,
                \Omega^2_\tX(-E_x) \tensor L'_x ) = 0.$$
 Now by Leray spectral sequence for $\pi_x:\tX \to X$, we
see that the above $H^1$ is vanishes if:
\begin{description}
    \item[{\rm(\thepropose.7)}]  $H^1(X , (\pi_x)_*(
                \Omega^2_\tX(-E_x) \tensor L'_x)) = 0$, and
    \item[{\rm(\thepropose.8)}] $R^1(\pi_x)_*(\Omega^2_\tX(-E_x)
                \tensor L'_x) = 0$.
 \end{description}
Note that by the projection formula,
        $$(\pi_x)_*(\Omega^2_\tX(-E_x)
        \tensor L'_x) = (\pi_x)_*(\Omega^2_\tX(-2E_x)) \tensor L.$$
 Hence vanishing of \specialref{sublemma1}{7} above is
implied by hypothesis \specialref{base point and smooth}{1}. So
we have to check that $R^1$ is also zero. By projection formula, it
suffices to prove the following:
\begin{equation}\label{whats in a name}
        R^1(\pi_x)_*(\Omega^2_\tX(-2E_x)) = 0.
\end{equation}
Since the fibres of $\pi_x$ have dimension $\leq 2$, $R^2$ satisfies base
change.
Since the $R^2$ is supported on the point $x$, we can compute the
$R^2$ by restricting to $E_x$. Providing we show that $R^2$ is 0, we
get that $R^1$ satisfies base change. Then $R^1$ is also supported on
the point $x$ and can be computed by restricting to $E_x$.
On $E_x$ we have the fundamental exact
sequence
        $$ 0 \to \O_{E_x}(1) \to \Omega^1_\tX\big|_{E_x} \to
                \Omega^1_{E_x} \to 0.$$
An easy calculation shows that this exact sequence splits. So that we have
        $$\Omega^1_\tX\big|_{E_x} = \Omega^1_{E_x} \oplus
                \O_{E_x}(1).$$
So that by taking exteriors, we have
        $$\Omega^2_\tX\big|_{E_x} = \Omega^2_{E_x} \oplus
                \Omega^1_{E_x}(1).$$
 Hence we have
        $$\Omega^2_\tX(-2E_x)\big|_{E_x} =
                \Omega^2_\tX\big|_{E_x}\tensor \O_{E_x}(2) =
                        \Omega^2_{E_x}(2) \oplus \Omega^1_{E_x}(3).$$
Since $E_x \isom \P^2$, we see easily that $H^1$ and $H^2$ of the
latter sheaf vanish.
 \end{proof}

\begin{sublem}
    Hypothesis \specialref{base point and smooth}{2} implies
\specialref{proof of 3.2.3}{6}
 \end{sublem}
\begin{proof}
        We want to prove that
        $$H^2( \tX , \Omega^1_{(\X,\Y)/\tX}(D)\big|_{\X_s}) = 0.  $$
 Firstly we recall that $\Omega^1_{(\X,\Y)/\tX} \big|_{\X_s} =
M(L'_x, W)^* \tensor L'^{-1}_x$ (this is just the formula \ref{from
Kapils paper}).
 Further $\O_\X(D)\big|_{\tX} = \O_\tX(E_x)$. So that we
have to show that
        $$ H^2( \tX , M(L'_x,W)^* \tensor \O_\tX(E_x) \tensor
                L'^{-1}_x ) = 0.$$

By Serre duality, we have to prove the following vanishing:
        $$H^1(\tX , M(L'_x,W) \tensor \O_\tX(-E_x) \tensor K_\tX
                \tensor L'_x ) = 0.$$

    This is done as before by a Leray spectral sequence argument.  We
observe that $ H^1( X , (\pi_x)_*( M(L'_x,W) \tensor \O_\tX(-E_x)
\tensor K_\tX \tensor L'_x)) = 0 $ by the hypothesis \specialref{base
point and smooth}{2}. Thus it suffices to prove that
        $$R^1(\pi_x)_*(M(L'_x,W) \tensor \O_\tX(-E_x) \tensor K_\tX
                \tensor L'_x) = 0.$$
Using the projection formula and the fact that $K_\tX = \pi_x^*(K_X)
\tensor \O_\tX(2E_x)$, we see that it suffices to prove that
        $$R^1(\pi_x)_*(M(L'_x,W)) = 0.$$

As before we will show vanishing of $H^1$ and $H^2$ after
after restricting to $E_x$. This will imply the vanishing of $R^1$ and
$R^2$ as before. Since we have the exact sequence
        $$ 0 \to M(L'_x , W) \to W\tensor \O_\tX \to L'_x \to 0, $$
restricting this exact sequence to $E_x$, by the universal property
of the tautological bundle $\O_{E_x}(1)$, we see that we have:
        $$M(L'_x,W)\big|_{E_x} = \Omega^1_{E_x}(1) \oplus V' \tensor
                \O_{E_x},$$
for some subspace $V' \subset W$. Now since $E_x\isom\P^2$
one notes that the $H^1$ and $H^2$ of the latter sheaf vanish.
 \end{proof}

Finally,  to complete the proof of Proposition \ref{base point
and smooth}, it remains to prove \specialref{base point and smooth}{4}:

\begin{lemma}\label{proof of 3.2.4}
    $$H^1(\X_s, \Omega^1_{(\X,\Y)}(D)\big|_{D}\big|_{\X_s}) =0.$$
 \end{lemma}

\begin{proof}
        Observe that as $D=E_x \times S$, we have $D\big|_{\X_s} =E_x$
and hence $D\big|_{E_x} = \O_{E_x}(-1)$. We have the
divisor $Z = E_x\times S\cap \Y = D\cap \Y$ in $D$.
 Note that $D$ is just a product of $\P^2$
and a projective space. Thus the vanishing to be proved is reduced to
a vanishing result on $\P^2$. Further note that by definition $Z =
\P_{E_x}(M(L'_x,W)^*\big|_{E_x})$. We have a natural commutative
diagram
        $$\let\scty=\scriptstyle
        \Matrix{
         &  & \scty{0} &&\scty{0} & & \scty{0} &  & \cr
         &  & \scty{\da} &&\scty{\da} & &\scty{\da} &  & \cr
        \scty{0} & \scty{\to} &\scty{\O_D(-Z)} &\scty{\to}&\scty{\O_D} &
\scty{\to} &\scty{\O_Z} &\scty{\to} & \scty{0}\cr
         &  & \scty{\da} && \scty{\da} & & \scty{\da} &  & \cr
\scty{0} & \scty{\to} & \scty{\Omega^1_{(\X,\Y)}(D)\big|_D} &\scty{\to}&
\scty{\Omega^1_\X(D)\big|_D}
                                &\scty{\to} &\scty{\Omega^1_\Y(D)\big|_D} &
 \scty{\to} & \scty{0}\cr
         &  & \scty{\da} && \scty{\da} & & \scty{\da} &  & \cr
\scty{0} & \scty{\to} & \scty{\Omega^1_{(D,Z)}\tensor \O_{E_x}(-1)} &
\scty{\to}&
   \scty{\Omega^1_D\tensor
        \O_{E_x}(-1)} &\scty{\to} &\scty{\Omega^1_Z\tensor \O_{E_x}(-1)} &
\scty{\to} & \scty{0}\cr
         &  & \scty{\da} && \scty{\da} & & \scty{\da} &  & \cr
         &  &\scty{0} && \scty{0} & & \scty{0} &  & \cr
                }
        $$

Further restricting the right column to $\X_s = \tX\times \{s\}$ we see that
the middle term
        $H^1(\tX , \Omega^1_{(\X,\Y)}(D)\big|_D\big|_{\X_s}) = 0$
 if

\begin{description}
    \item[{\rm(\thepropose.9}] $H^1(E_x , \O_D(-Z)\big|_{\X_s}) =0$, and
    \item[{\rm(\thepropose.10)}] $H^1( E_x , \Omega^1_{(D,Z)}\big|_{E_x}
 \tensor
                \O_{E_x}(-1)) = 0$.
 \end{description}

Note that as $E_x \isom \P^2$, and as $\O_D(-Z)\big|_{E_x}$ is a line
bundle on $E_x$, its $H^1$ is trivially zero. So we have to prove
\specialref{proof of 3.2.4}{10}. To do this we proceed as follows: We
have the exact sequence (from the definition)
        $$\Matrix{
        0 &\to& \Omega^1_{(D,Z)}& \to &\Omega^1_D &\to& \Omega^1_Z& \to& 0\cr
          &    &                &     & \Vert     &    &          & &\cr
          &    &  & & \Omega^1_{E_x}\oplus \Omega^1_S &    &      & &\cr
                }
        $$
   So restricting to $E_x$, as $\tilde{s}$ cuts out a smooth divisor
on $\tX$, $Z=E_x\cap \Y$ is a $\P^1$-bundle over $S$
in a Zariski neighbourhood of $\tilde s$.
 Then we have the exact sequence
$$
0 \to \Omega^1_{(D,Z)}\big|_{E_x} \to \Omega^1_D\big|_{E_x} \to
                \Omega^1_Z \big|_{E_x} \to 0.$$
Now $\Omega^1_D\big|_{E_x}=\Omega^1_{E_x}\oplus
\Omega^1_{S,s}\tensor\O_{E_x}$, and
similarly we have
        $\Omega^1_Z\big|_{E_x} = \Omega^1_{\P^1} \oplus
                \Omega^1_{S,s}\tensor\O_{E_x}$.
This gives an exact sequence
        $$\displaylines{
0  \to (\Omega^1_{(\P^2,\P^1)} \oplus \Omega^1_{S,s}\tensor\O_{E_x})
\tensor \O_{E_x}(-1)  \to
  \Omega^1_{\P^2}(-1) \oplus
                \Omega^1_{S,s}\tensor \O_{E_x}(-1)
        \to  \hfill\cr
        \hfill{}\qquad\to \Omega^1_{\P^1}(-1) \oplus
                \Omega^1_{S,s} \tensor \O_{E_x}(-1)\to 0.\hfill\cr
        }
        $$
Then taking long exact cohomology sequence we have:
$$\displaylines{
H^0(\P^1, \Omega^1_{\P^1}(-1)) \oplus H^0(\P^2,\Omega^1_{S,s}
 \tensor \O_{E_x}(-1)) \to \hfill\cr
        \hfill{}\qquad \to H^1( (\Omega^1_{(\P^2,\P^1)}
                \oplus \Omega^1_{S,s}\tensor\O_{E_x})\tensor
                \O_{E_x}(-1)) \to \hfill\cr
        \hfill{} \qquad\qquad \to H^1(\Omega^1_{\P^2}(-1))
 \oplus H^1(\Omega^1_{S,s}\tensor\O_{E_x}(-1) ).\hfill\cr
        }
$$

And the terms on the extreme vanish so that the middle term vanishes.
Hence we are done. This finishes the proof of Lemma \ref{proof
of 3.2.4}
 \end{proof}

        For the next proposition, we need some more notation. Let
        $$ L''_x = (\pi_x)^*(L)\tensor\O_\tX(-2E_x).$$
 Observe that image, under $\pi_x$ of smooth sections of
$H^0(\tX,L''_x)$, are sections in $H^0(X,L)$ which have a single
ordinary double point at $x$. We now want to reduce the condition
$\vnl(\tX,L''_x,H^0(\tX,L''_x),s)$ for smooth $s\in
S=\P(H^0(\tX,L''_x)^*)$, to a more manageable form. This is done by
the following.

\begin{propose}\label{base point and ODP}
    Assume that $L$ is such that $L''_x$ is very ample on $\tX$.
Let $s\in W=H^0(\tX, L''_x)$ be a section which cuts out a smooth
divisor  on $\tX$. Assume that the following assertions are valid:
        \begin{description}
            \item[{\rm(\thepropose.1)}] $H^1(X, L\tensor m_x^3) = 0$,
                where $m_x$ is the maximal ideal of $x\in X$.
            \item[{\rm(\thepropose.2)}] $H^1( X ,
                (\pi_x)_*(\Omega^2_{\tX_x}(-2E_x)) \tensor L) = 0$, and
            \item[{\rm(\thepropose.3)}] $H^1( X , (\pi_x)_*( M(L_x'') )
                \tensor K_X \tensor L ) = 0$.
         \end{description}
Then $\vnl(\tX,L''_x,H^0(\tX,L''_x),s)$ is valid.
 \end{propose}

\begin{proof}
    Since the linear system $W$ is base point free on $\tX$, and as
the section $s\in W$ cuts out a smooth divisor on $\tX$, we see that
Proposition \ref{base free and smooth} can be applied. By this
proposition, it suffices to prove the following assertions:

\begin{description}
    \item[{\rm(\thepropose.4)}] $H^1(\tX, \Omega^2_\tX \tensor L''_x ) = 0$,
and

    \item[{\rm(\thepropose.5)}] $H^1( \tX , M(L''_x,W) \tensor K_\tX \tensor
                L''_x ) =0$.
 \end{description}

Thus we have to prove that the hypothesis of the proposition ensure
the vanishing \specialref{base point and ODP}{4} and
\specialref{base point and ODP}{5}. The argument is similar to the
one earlier. Firstly, let us observe that the implication
\specialref{base point and ODP}{2} implies \specialref{base
point and ODP}{4} follows from  \ref{sublemma1}. By a Leray
spectral sequence argument, we are reduced to proving
        $$R^1(\pi_x)_*(\Omega^2_\tX(-2E_x)) = 0.$$
 But this has been proved during the proof of
Sublemma\ref{sublemma1} as equation (\ref{whats in a name}).
 Thus we are done in this case. So it
remains to prove that \specialref{base point and ODP}{1} and
\specialref{base point and ODP}{3} together imply
\specialref{base point and ODP}{5}.

First we check that $R^1(\pi_x)_*(M(L''_x,W)) = 0$. This can be done
by restricting to $E_x$. On $E_x$ we have an exact sequence:
        $$ 0  \to M(L''_x,W)\big|_{E_x} \to W\tensor \O_{E_x} \to L''_x
                \big|_{E_x} \to 0.$$
So on noting that $W= H^0(\tX,L''_x)=H^0(X,L\tensor m_x^2)$, the
result now follows from \specialref{base point and ODP}{1}
and \specialref{base point and ODP}{3}.
 \end{proof}

\section{The Noether-Lefschetz Locus}\label{The Noether-Lefschetz Locus}
        For any smooth projective threefold $X$, and $L$ be a very
ample line bundle on $X$, we say smooth member $Y$ of $H^0(X,L)$ lies
in the ``Noether-Lefschetz locus of the linear system $H^0(X,L)$'' if
$\Pic(X)\to \Pic(Y)$ is not a surjection. More generally, if $Y$ is
normal we say $Y$ lies in the ``Noether-Lefschetz locus'' if $\Pic(X)
\to CH^1(Y)$ is not a surjection.

        Now we  deduce the global Noether-Lefschetz theorem from the
infinitesimal results.

\begin{propose}\label{codim one}
    Suppose $X$ is a smooth projective threefold over $\C$ and $L$ is
a very ample line bundle over $X$. Assume that $\inl(X,L,H^0(X,L),s)$
is valid for all $s\in H^0(X,L)$ which cut out a smooth divisor on
$X$. Then the Noether-Lefschetz locus for the linear system
$H^0(X,L)$ has codimension $\geq 1$ in $\P(H^0(X,L)^*)$.
 \end{propose}

\def\bK{\bar{K}}

\begin{proof}
    We argue as in \cite{Mohan-Srinivas}. We can assume that the $X,
L$ are defined over a finitely generated field $K/\Q$. Let $\bK$
denote the algebraic closure of $K$.  Let $\bK(\eta)$ denote the
rational function field in $\dim H^0(X,L)$ variables over $\bK$.  Then
this defines a point $\eta \into S=\P(H^0(X,L)^*)$, which we will call
the $\bK$-generic point of $S$.  Let $\etab \in S=\P(H^0(X,L)^*)$ be
the corresponding $\bK$-geometric generic point of $S$.  Then we claim
that $\Pic(X_\etab) \isom \Pic(Y_\etab)$.

        The claim is proved as follows. Suppose that it is not an
isomorphism, then in particular the map is not surjective. Let
$\alpha\in\Pic(Y_\etab)$ be any line bundle not in the image of the
map. Then this cycle is defined over a finite field extension of the
function field $\bK(\eta)$. Then there exists an \'etale open set $U
\to S$ such that $\alpha$ ``spreads'' out to a line bundle on $\Y_U$.
By replacing  $U$ by a smaller non-empty open subset if necessary, we
can further assume that $\alpha$ is a nontrivial element of
$\Pic(\Y_U)$. Now we restrict $\alpha$ to the fibres of $\Y_U \to U$.
By shrinking $U$ further we can assume that image of $U$ in $S$ is an
open set which does not meet the discriminant locus of the linear
system $H^0(X,L)$ in $S$. Thus we can now apply the Infinitesimal
Noether-Lefschetz to the fibres over $U$. On the fibres of $\Y_U \to
U$, as $\inl(X,L,H^0(X,L),s)$ is valid, we see that the restriction of
$\alpha$ to the fibres over $u\in U$ is zero in $\Pic(Y_u)$ because
infinitesimally there are no extra cycles (note that we have used the
fact that we can identify $\O_U/m_u^2 \isom \O_S/m_s^2$, where $s$ is
the image of $u$ under the map $U\to S$). Then by the semi-continuity
theorem we see that the line bundle $\alpha$ on $\Y_U$ must be
trivial. This proves the claim.

Now it remains to prove that the Noether-Lefschetz locus has
codimension $\geq 1$. We proceed as in \cite{Mohan-Srinivas}. Let
$U$ be the subset of $S$ obtained by removing all the divisors of
$S$, which are defined over $\bK$. Thus we have removed a
countable set of closed subvarieties. If $p\in U$ is a closed point
the map $ \Spec(\C(p)) \into U$ factors through the geometric generic
point $\etab$. Further as $\Y_\etab$ is the zero scheme of a section
of an ample line bundle on $X_\etab$, we see as $H^1(X_\etab,
\O_{X_\etab}) \isom H^1(Y_\etab,\O_{Y_\etab})$. We see that
$\Pic(Y_p)(\bK)\isom \Pic(Y_\etab)(\bK)$. This proves the claim.
  \end{proof}

        An identical argument, with divisors on $S$ which are not
contained in the support of the discriminant locus, we can prove the
following:

\begin{propose}\label{codim two in smooth}
    If {}$\inl(X,L,W,s)$ is valid for all $W$ of codimension one in
$H^0(X,L)$ and $s\in W$ smooth. Then ``smooth part'' of the
Noether-Lefschetz locus has codimension $\geq 2$ in $S$.
 \end{propose}

\begin{proof}
    Since the proof is almost identical to the one given earlier, we
will only indicate the important point. The idea is to work with
divisor of $D$ which is not contained in the discriminant locus. For
such a divisor, we prove by an argument identical to the one given
above that for the $\bK[A$-geometric generic point $\etab_D$ of $D$, we
have an isomorphism $\Pic(X_{\etab_D}) \isom \Pic(Y_{\etab_D})$.  This
is proved as before, except that one has to use $\inl$ for a
codimension one linear system. For which we use the following Lemma.
Assuming the Lemma for the moment, we  complete the proof as follows.
Consider the subset of $D$ obtained by removing all the codimension
one subvarieties of $D$ which are defined over $\bK$. Then the
one argues as in the proof of Proposition~\ref{codim one}.
 \end{proof}

\begin{lemma} \label{vnl and inl on the blowup}
    Let $(X,L,W,s)$ be such that $W\subset H^0(X,L)$ is linear system
of codimension one with a single base point at $x\in X$. Further
assume that $s$ cuts out a smooth divisor $Y_s$ in $X$. Then we have
        $$ \inl(\tX_x,L'_x,W,\tilde s) \implies \inl(X,L,W,s).$$
 \end{lemma}

\begin{proof}
        We use the following commutative diagrams:
        $$
        \Matrix{
\Pic(\tX_x)&& \mapright{f_2} && \Pic(\tilde{\Y_T}) \cr
 h_2 \searrow   &&      && \swarrow g_2      \cr
                && \Pic(\tY_s)  &&           \cr
 \mapup{\tau_1}   && \mapup{\tau_2}     &&  \mapup{\tau_3}     \cr
                &&    \Pic(Y_s)    &&     \cr
 h_1 \nearrow    &&      && \nwarrow g_1      \cr
\Pic(X)       && \mapright{f_1} &&   \Pic(\Y_T)    \cr
                }
        $$
And one also has the following commutative diagram:
$$
        \Matrix{
0 & \to & \Pic(X) & \mapright{\tau_1} & \Pic(\tX_x) &
                \mapright{\tau_1'}& \Z & \to &  0                \cr
  &     & \mapdown{h_1} &             & \mapdown{h_2}&
                                                    &\Big\Vert& & \cr
0 & \to & \Pic(Y_s) & \mapright{tau_2}& \Pic(\tY_s) &
                                \mapright{\tau_2'} & \Z & \to & 0 \cr
                }
$$
        As $\inl(\tX_x,L'_x,W,s)$ is valid, we have
$\image{(h_2)}=\image{(g_2)}$. We want to prove that
$\image{(h_1)}=\image{(g_1)}$. Since the lower triangle commutes, we
have the inclusion $\image{(h_1)} \subset \image{(g_1)}$. So we will
prove the reverse inclusion. Let $\alpha\in \Pic(\Y_T)$, as $\tau_2
g_1 = g_2\tau_3$ we have $\tau_2g_1(\alpha) = g_2\tau_3(\alpha)$. This
implies that there exists $\alpha'\in \Pic(\tX_x)$ such that
$h_2(\alpha') = g_2(\tau_3(\alpha)) = \tau_2(g_1(\alpha))$.

        Now chasing the image $g_1(\alpha)\in\Pic(Y_s)$ in the second
diagram, we see that  $\tau_2'(h_2(\alpha')) =
\tau_2'(\tau_2(g_1(\alpha))) =0$. Then $\tau_2'(h_2(\alpha')) =
\tau_1'(\alpha') = 0$. So that $\alpha' \in \image(\tau_1)$. So there
is an $\alpha'' \in \Pic(X)$ such that $\alpha' = \tau_1(\alpha'')$.

        Now going back to the previous diagram, we see that
        $$h_2(\alpha')=h_2(\tau_1(\alpha'')) = \tau_2(h_1(\alpha''))=
                \tau_2(g_1(\alpha)).$$
 Hence, $\tau_2(h_1(\alpha'')-g_1(\alpha)) =0$. But as $\tau_2$ is
injective, we see that $h_1(\alpha'')=g_1(\alpha)$.

        This proves Lemma~\ref{vnl and inl on the blowup}.
\end{proof}

Applying the same technique we can prove:

\begin{propose}\label{codim one in disc.}
    If for all $x\in X$, $\inl(\tX,L''_x,H^0(\tX,L''_x),s)$ is valid
for all smooth $s\in H^0(\tX,L''_x)$, then the intersection of the
Noether-Lefschetz locus of $H^0(X,L)$ with the discriminant locus of
$H^0(X,L)$ has codimension $\geq 2$ in $S=\P(H^0(X,L)^*)$.
 \end{propose}

\begin{proof}
    Firstly one notes that the hypothesis of the proposition together
with Proposition~\ref{codim one}, implies that the
Noether-Lefschetz locus of $H^0(\tX,L''_x)$ is a proper of
codimension $\geq 1$ in $\P(H^0(\tX,L''_x)^*)$.

Note that the image of the zero scheme of any smooth $s\in
H^0(\tX,L''_x)$ in $X$ is an element of the linear system $H^0(X,L)$
with a single ordinary double point at $x$. Then the intersection of
the Noether-Lefschetz locus with the discriminant locus is clearly of
codimension $\geq 2$. This proves the result.
 \end{proof}

\section{A General Noether-Lefschetz Theorem}\label{A General
Noether-Lefschetz Theorem}

In this section we are going to apply the general machinery of the
previous section to prove the following generalization of the
classical Noether-Lefschetz theorem.

\begin{thm}\label{main theorem}
    Let $X/\C$ be a smooth projective threefold, $L$ a very ample
line bundle on $X$. Then there exits an integer
$n_0(X,L) > 0$  such that for all $n \geq n_0(X,L)$, the
Noether-Lefschetz locus of the linear system $H^0(X,L^n)$ has
codimension $\geq 2$.
 \end{thm}

\begin{proof}
    The general machinery for the proof of this result was set up in
the previous section. By the Propositions \ref{base free and
smooth}, \ref{base point and smooth} and \ref{base point and ODP},
and Propositions \ref{codim one}, \ref{codim two in smooth} and
\ref{codim one in disc.} we are reduced to proving the following
assertion:

There exists and $n_0(X,L)$ such that for all $n\geq n_0$ the
following are valid:

\begin{enumerate}
    \item $H^1(X,\Omega^2_X\tensor L^n) =0$,
    \item for every base point free linear system $W\subset
        H^0(X,L^n)$ of codimension at most one, we have
                $$H^1(X,M(L^n,W) \tensor K_X \tensor L^n ) =0.$$
    \item For all points $x\in X$, for all $n \geq n_0$ we have
                $$H^1( X , (\pi_x)_*( \Omega^2_{\tX_x}(-2E_x)) \tensor
                        L^n ) = 0, $$
    \item for any codimension one linear system $W\subset H^0(X,L^n)$
                with a single base point at $x\in X$ the following
                holds for all $n \geq n_0$:
                $$H^1( X , (\pi_x)_*(M(L'_x,W))  \tensor K_X
                        \tensor L^n) =0.$$
     \item For all $x\in X$ and for all $n\geq n_0$, we have
                $$H^1(X , m^3_x \tensor L^n ) =0,$$
     \item
                $$H^1(X, (\pi_x)_*(\Omega^2_{\tX_x}(-2E_x)) \tensor
                L^n ) = 0;$$
     \item for $W= H^0(\tX_x,L''_x)$, we have for all
                $x\in X$ and for all $n\geq n_0$:
                $$H^1( X , (\pi_x)_*(M(L''_x,W)) \tensor K_X \tensor
                L^n ) = 0.$$
 \end{enumerate}

        By what has transpired so far, these assertions, by the
results of Section~\ref{The Noether-Lefschetz machine}, imply the
infinitesimal Noether-Lefschetz results \ie Propositions \ref{base
free and smooth}, \ref{base point and smooth} and \ref{base point
and ODP}, and finally from these local results by
Propositions \ref{codim one}, \ref{codim two in smooth} and
\ref{codim one in disc.}, we obtain the above global result. This
entire section is devoted to the proof of these seven assertions.
Firstly, note that of these seven assertions, the assertions (1), (3),
(5) and (6) follow immediately from the Semi-continuity Theorem,
Serre's vanishing theorem and from simple Noetherian induction. This
will be left to the reader. The remaining assertions are more
difficult and require more elaborate arguments. We will now prove the
assertions (2), (4) and (7).

 \end{proof}

\begin{propose}
    There exists a positive integer $n_0$ depending only on $X,L$ such
that for any base point free linear system of codimension at most
one, we have
        $$H^1( X , M(L^n,W)\tensor K_X \tensor L^n ) = 0.$$
 \end{propose}

\begin{proof}
        The idea is to use a method of Mark Green (see
\cite{MarkGreen}) and reduce the result to a regularity computation.

    Let $V = H^0(X,L)$, by hypothesis, $X$ embeds in $\P:=\P(V)$.
Let $j:X \into \P$ be the embedding given by the linear system. Write
$\I_X$ for the ideal sheaf of $X$ in $\P$.
Then note that $L = j^*(\O_\P(1))$. We have a surjection:
        $$H^0( \P, \O_\P(n)) \to H^0(X,L^n) \to 0.$$
 We write $F=H^0(X,L^n), F'=H^0(\P,\O_\P(n))$. Then the subspace
$W\subset F$ gives us a subspace $W'\subset F'$. If we choose $n_0 >
\reg(\I_X)$ then $W'$ is base point free on $\P$. Moreover, one has
$\codim_F'(W') = \codim_F(W)$. Then we have corresponding to the
triple $(\P,\O_\P(n),W')$ a vector bundle $M':=M(\O_\P(n),W')$ which
is defined by the usual evaluation sequence:
        $$ 0 \to M' \to F'\tensor \O_\P \to \O_\P(n) \to 0.$$

         Write $V' = \ker(W' \to W)$. Now we have the following
commutative diagram of locally free sheaves on $X$:
 $$
\Matrix{
  &     &       0    &     &    0          &     &      &            &   \cr
  &     &      \da   &     &  \downarrow   &     &      &            &   \cr
  &    &V'\tensor\O_X&  =  & V'\tensor\O_X &     &      &            &   \cr
  &     &   \da      &     & \downarrow    &     &      &            &   \cr
0 & \la & j^*M'      & \la & W'\tensor\O_X & \la &    F &    \la     & 0 \cr
  &     & \downarrow &     & \downarrow    &     & \parallel &       &   \cr
0 & \la & M          & \la & W\tensor\O_X  & \la &    F &        \la & 0 \cr
  &     & \downarrow &     & \downarrow    &     &      &            &   \cr
  &     &    0       &     &    0          &     &      &            &   \cr
 }
$$

 Now any splitting of the middle column gives a (non-canonical)
splitting of the left column. Thus we see that one has a
non-canonical splitting:

         $$  j^*(M') = M(L^n,W) \oplus V'\tensor \O_X.$$

        Now we need the following lemma due to Mark Green,
see \cite{MarkGreen}.

\begin{lemma}
    Let $W$ be any base point free linear system $W\subset
H^0(\P,\O_\P(n))$. Then for all $i \geq 1$, and for all $k+i \geq
\codim(W)+1$, we have
        $$H^i( \P , M' \tensor \O_\P(k)) =0,$$
 in other words, $M'$ is $(\codim(W)+1)$-regular.
 \end{lemma}

 Now by projection formula we have $j_*j^* M' =j_*M(L^n,W) \oplus V'
\tensor j_*(\O_X)$. Thus as $j_*(M(L^n,W))$ is a direct summand of $M'
\tensor j_*(\O_X)$, and as the regularity of $M'$ is bounded
independent of $n$, and as the regularity of $j_*(\O_X)$ is a fixed
constant, it follows that the regularity of $j_*M(L^n,W)$ is bounded
independent of $n$. Thus if we assume that $n >
(\codim(W)+1+\reg(K_X))$, then $H^1(X, M(L^n,W)\tensor K_X \tensor L^n
) =0$. This proves the proposition.

 \end{proof}

        A slight modification of the same technique gives the
following result.
\begin{propose}
    For any codimension one linear system $W\subset H^0(X,L^n)$ with a
single base point $x$, there is an $n_0(X,L)$ such that for all
$n\geq n_0$, we have:
        $$H^1( X , (\pi_x)_*(M(L'_x,W)) \tensor K_X \tensor L^n )
                =0.$$
 \end{propose}

\begin{proof}
    We have to modify the argument we gave for the previous proposition.
We need the following variant of the previous lemma.

        For $x\in \P$ write $\tilde\P$ for the blowup of $\P$ at $x$,
and write $\pi_x$ for the blowup morphism. Write
$W=H^0(\P,\O_\P(n)\tensor m_{x,\P})$, where $m_{x,\P}$ is the maximal
ideal of $x$ in $\P$. $W$ gives a natural linear system on $\tilde\P$.
Further, let $M'$ be the vector bundle defined in
Section~\ref{Generalities and notations}, on $\tilde\P$ corresponding
to the linear system $W$. Assume that $n \geq 2$ so that this linear
system is base point free. We have

\begin{lemma}
With the above notations, the coherent sheaf $(\pi_x)_*(M')$ on $\P$
is $1$-regular, \ie for all $i\geq 1$ and for all $k+i \geq 1$ we have

        $$H^i(\P, (\pi_x)_*(M') \tensor \O_\P(k)) = 0.$$
 \end{lemma}

\begin{proof}
    On $\P$ we have the exact sequence
        $$ 0 \to (\pi_x)_*(M') \to W\tensor \O_\P \to \O_\P(n) \tensor
                m_{x,\P} \to 0.$$
 The result now follows from the regularity computation for the sheaf
$m_{x,\P}$.
 \end{proof}

\long\def\comment#1\endcomment{}
\comment
        Note that we have assumed here that $n$ is large enough, so
that the $W$ spans $\O_\P(n) \tensor m_{x,\P}$. Now twisting by
$\O_\P(k-i)$ we see that we have
        $$\to H^{i-1}(\O(n+k-i) \tensor m_x ) \to H^i( (\pi_x)_*(M')
                \tensor \O(k-i)) \to H^i(\O(k-i)) \tensor W\to$$
         Now for $i=1$, $H^0(\O(k-1)) \tensor W \to
H^0(\O(n+k-1)\tensor m_x)$ is a surjection if $k\geq i+1$. And for
$1\leq i \leq \dim\P$, we have $H^i(\P,\O(k-i) =0$ and for $i=\dim\P$,
$H^i(\O_\P(k-i)) = 0$ if $k-i > -\dim\P-1$ \ie if $k > -1$. So one has
to prove that the first term on the left in the exact sequence above
is zero. This is done as follows: we have the exact sequence
        $$ 0 \to m_x \to \O_\P \to \C_x \to 0,$$
 where $\C_x$ is the skyscraper sheaf supported only at $x$. Then
twisting by $\O(n+k-i)$ we get:
        $$ 0 \to \O(n+k-i)\tensor m_x \to \O_\P(n+k-i) \to
                \C_x\tensor\O(n+k-i)  \to 0.$$
         Note the sheaf on the extreme right is still supported at a
point and hence has no cohomology except possibly $H^0$. Hence for
$i\geq 1$ we have $H^i(m_x\tensor\O(n+k-i)) \isom H^i(\O(n+k-i))$. So
that for $i=\dim\P$, we have $H^i(\O(n+k-i))=0$ if $n+k-i > -\dim\P-1$
\ie if $n+k > -1 $. Thus  finally if $k > -1 $ and $k > -1 -n $ then
for $i=\dim\P$, $H^i( (\pi_x)_*(M') \tensor \O(k) ) = 0 $ if $k+i \geq
1$. Thus we have proved that $(\pi_x)_*(M')$ is $1$-regular.
 \end{proof}
\endcomment

        Finally, we can now finish the proof of the proposition. We
have an exact sequence
        $$ 0 \to \pi_*M' \to \O_\P \tensor H^0(\O(n)\tensor m_{x,\P})
                \to \O(n) \tensor m_{x,\P} \to 0.$$

         Write $K_x = \ker(L^n \tensor m_{x,\P} \tensor \O_X \to L^n
\tensor m_{x,X})$. Note that $K_x$ is a finite length sheaf, supported
at the point $x$. We assume that $n$ is large enough to ensure that
there is a surjection
        $$
        \Matrix{
H^0(\P, \O(n)\tensor m_{x,\P}) & \to & H^0(X, L^n\tensor m_{x,X})& \to& 0\cr
          \Vert                &     &       \Vert               &    & \cr
           F'                  &     &       F                   &    & \cr
                }
        $$
         Write $W'=\ker(F' \to F)$. Then there is a surjection
$W'\tensor \O_X \to K_x$. To see this note that $K_x$ is a sheaf
supported at a point, further note that we can assume, by choosing $n$
sufficiently large so that $H^1(X, L^n \tensor I_X ) = 0$, that $F'$
surjects on to $H^0(X, L^n \tensor m_{x,\P} \tensor \O_X)$. This being
done, we have the required surjection. Note that this choice of $n$
depends only on $X,L$. We have then the following commutative diagram
on $X$

$$\let\sc=\scriptstyle
\Matrix{
  & & \sc{0} & &                \sc{0} &  & \sc{0} &  & \cr
  & & \sc{\da} & &              \sc{\da} &  & \sc{\da} &  & \cr
 \sc{0} & \sc{\to}& \sc{H} &\sc{\to}&\sc{W'\tensor\O_X} & \sc{\to} & \sc{K_x} &
   \sc{\to} & \sc{0}\cr
  & & \sc{\da} & &              \sc{\da} &  & \sc{\da} &  & \cr
 \sc{0} & \sc{\to}& \sc{j^*M'/ \Tor^{\O_X}_1( \C_x, \O_X)} &\sc{\to}&
                \sc{F'\tensor\O_X} &\sc{\to} & \sc{L^n \tensor m_{x,\P} \tensor
                        \O_X} & \sc{\to} & \sc{0}\cr
 & & \sc{\da} & &               \sc{\da} &  & \sc{\da} &  & \cr
 \sc{0} & \sc{\to}& \sc{M} & \sc{\to}&\sc{F\tensor\O_X} & \sc{\to} &
\sc{m_{x,X}
   \tensor L^n} & \sc{\to} & \sc{0}\cr
 & & \sc{\da} & &               \sc{\da} &  & \sc{\da} &  & \cr
 & & \sc{0} & & \sc{0} &  & \sc{0} &  & \cr
        }
$$

         Note that $\Tor^{\O_X}_1(\C_x,\O_x)$ is supported at $x$ and
consequently is of finite length. Further we note that $K_x$ is
also a sheaf of finite length.

        Now we want to compute the regularity of $M$ in terms of
$j^*M'$. On one hand we know that the regularity of $j^*M'$ can be
computed in terms of regularity of $M'$ on $\P$. By the previous lemma,
we know that the regularity of $M'$ on $\P$ is independent of $n$. And
so that regularity of $j^*M'$ depends only on regularity of $M'$ and
on regularity of $j_*\O_X$. Thus we can choose $n$ to be larger than
the sum of these two numbers.

        We note that as the $\Tor$ term is of finite length, the
regularity of the sheaf $j^*(M')$ is the same as the regularity of the
sheaf $j^*(M')/\Tor^{\O_X}_1(\C_x,\O_X)$ (this is essentially because
the $\Tor$ has no higher cohomologies, being supported on a single
point). Similarly, as $K_x$ has finite length, the top row of the
diagram says that $F$ and $\O_X$ have the same regularity. Thus we see
that there is an $n_0$ depending only on $X,L$ such that $M$ is $n_0$
regular.  Then the vanishing which is required follows.
 \end{proof}

        Lastly, the assertion (7) is proved exactly as above. The
point to be observed is the following lemma

\begin{lemma}
    Let $W=H^0(\P,\O(n) \tensor m_x^2 )$, $x\in\P$, define the
coherent sheaf $M''$ on $\P$ by the exact sequence:
        $$ 0 \to M'' \to H^0(\P,\O(n)\tensor m_x^2) \tensor \O_\P \to
                \O(n) \tensor m_x^2 \to 0.$$
 Then $M''$ is $2$-regular on $\P$.
 \end{lemma}

\section{Application to codimension two cycles}\label{Applications
to Codimension two cycles}
        In this section we give an application of our general
Noether-Lefschetz Theorem to codimension two cycles. The application
which we have, is related to a conjecture of Madhav Nori (see
\cite{Madhav}, \cite{KP}).  Let $X$ be a smooth projective variety
and $L$ a very ample line bundle on $X$. Let $S=\P(H^0(X,L^n)^*)$. Let
$K$ be the function field of $S$, and let ${\bar K}$ be its algebraic
closure. Let $\X=X\times S$ and let $\Y=\P(M(L^n)^*)\into \X$ be the
universal hypersurface corresponding to sections of $L^n$. The natural
inclusion $\Spec(K) \into S$ gives rise to a fibre square
        $$\Matrix{
                 \Y_K & \into & \X_K \cr
                  \downarrow & & \downarrow \cr
                  \Y & \into & \X \cr
                 }
        $$
The inclusion $\Y_K \into \X_K$ then gives rise to the restriction
map on rational equivalence classes of cycles in codimension $i$:
                $$ CH^i(\X_K) \to CH^i(\Y_K).$$
        Nori's conjecture in this set up is the following:

\begin{conj}[M. V. Nori]\label{Madhav's conjecture}
 If $n$ is sufficiently large then the natural map
        $$CH^i(\X_{\bar{K}})\tensor\Q \to CH^i(\Y_{{\bar K}})\tensor\Q$$
 is an isomorphism for $i < \dim(\Y_K) $ and an inclusion for
$i=\dim(\Y_K) $.
 \end{conj}

        Observe that the inclusion of fields $\C \into K$ gives rise
to a morphism $\X_K \to X_{\C}$ and hence a pull back morphism
$CH^i(X_{\C}) \to CH^i(\X_K)$ on cycle classes. Hence by composition
we get an homomorphism $p_X^* : CH^i(X_{\C}) \to CH^i(\Y_K)$.

        We can now prove the following weaker statement:

\begin{thm}\label{codim two cycles}
        Let $X$ be a smooth projective threefold, $L$ an ample line
bundle on $X$.  If $n$ is sufficiently large then the natural
map (given above):
        $$CH^2(X_{\C}) \to CH^2(\Y_K)$$
 is an inclusion.
 \end{thm}

         Before we begin the proof of the theorem, let us note one
immediate consequence of Theorem~\ref{codim two cycles}.

\begin{cor}
    Let $X$ be a smooth projective threefold, $L$ a very ample line
bundle on $X$. Then we have
        $$ CH^2(X_{\C})\tensor\Q \into CH^2(\Y_{\bar K})\tensor\Q $$
 \end{cor}

        This is immediate from the fact that the passage from $K$ to
${\bar K}$, annihilates only the torsion in Chow groups.

\begin{proof}
        Let $z\in CH^2(X_{\C})$ be  a cycle such that $p_X^*(z) = 0$
in $CH^2(\Y_K)$.  Then we note that there exists a divisor $D\into S$
such that $z$ is the image of a cycle in $CH^1(\Y_D)$, where $\Y_D$ is
the pull back of $\Y \to S$ to $D$ (via the inclusion $D\into S$).
We can find a  general pencil $\P^1 \into S$, which meets $D$
transversally. Let $F$ denote the function
field of the $\P^1$ which parametrises the pencil.  We have
a commutative diagram in which the rows are complexes:
        $$\Matrix{
                & &  CH^2(X_{\C}) && \cr
                & & \mapdown{\tau} & & \cr
     CH^1(\Y_D) & \mapright{f} & CH^2(\Y) & \mapright{h} & CH^2(\Y_K) \cr
     \mapdown{g} &     & \mapdown{g'} & &  \cr
     CH^1(\Y_{D\cap\P^1}) & \mapright{f'} & CH^2(\Y_{\P^1}) &
                \mapright{h'} & CH^2(\Y_F) \cr
                }
        $$
        Since $p_X^*(z) = h(\tau(z)) = 0$, there is a $\psi\in
CH^1(\Y_D)$ such that $f(\psi) = \tau(z)$. Then as the square commutes
and as the bottom row is a complex, we see that the image of $z$ in
$CH^2(\Y_F)$ is also zero. Thus we are reduced to proving the
following stronger assertion.
 \end{proof}

\begin{thm}\label{reduction to pencil}
        Suppose $n$ is sufficiently large, and let
$S=\P(H^0(X,L^n)^*)$, and let $\Y_P \into X\times \P^1$ be a general
Lefschetz pencil in $S$.  Let $F$ be the function field of the $\P^1$
corresponding to this pencil.  Then the natural morphism $CH^2(X_{\C})
\to CH^2(\Y_F)$ is injective.
 \end{thm}

\begin{proof}
        Since $n$ is sufficiently large, by Theorem~\ref{main
theorem} the Noether-Lefschetz locus is of codimension at least two in
$S$. Therefore we can assume that for every closed point $t\in\P^1$,
we have a surjection $\Pic(X) \to CH^1(\Y_t)$.  Moreover by  Bertini's
Theorem we can also assume that the base locus of such a pencil is an
irreducible, smooth projective curve $C$.

        Let $E\into \Y_P$ be the exceptional divisor. Observe that
$E=C\times \P^1$. Then it is easy to check that we have an isomorphism
$CH^2(\Y_P) \isom CH^2(X_{\C}) \oplus CH^1(C)$ (see \cite{Fulton}).
Further, we also have the following exact sequence:
        $$
        \oplus_{t\in\P^1}CH^1(\Y_t)  \to  CH^2(X)\oplus
                        CH^1(C)  \to  CH^2(\Y_F)
        $$
 and where the direct sum  extends over all the closed points of
$\P^1$.

        Let $z\in CH^2(X_\C)$ be  such that $p_X^*(z)=0$ in
$CH^2(\Y_F)$. Then $p_X^*(z)$ is in the image of the map
$\oplus_{t\in\P^1}CH^1(\Y_t) \to CH^2(\Y_P)$. And hence there are
finitely many  closed points $t_1,t_2\ldots,t_m\in \P^1$, and cycles
$s_i \in CH^1(\Y_{t_i})$ such that $p_X^*(z)=\sum_{i}s_i$.  By the
choice of our pencil, we see that there are cycles $w_i\in CH^1(X)$
such that $s_i= p_X^*(w_i).\Y_{t_i}$. But any two fibres of $\Y_P \to
\P^1$ are rationally equivalent. Hence $p_X^*(z) =
p_X^*(\sum_{i}w_i).\Y_t$.  Writing $w=\sum_{i}w_i$, we can rewrite the
last equality as $p_X^*(z) = p_X^*(w).\Y_t$.

        So to prove that $z =0$, by the projection formula, it
suffices to prove that $w=0$. Assume, if possible, that $w$ is not
zero. Then we have on intersecting with $E$, $p_X^*(z).E =
p_X^*(w).\Y_t.E$. By the projection formula, $(p_X)_*(p_X^*(z).E) =
z.(p_X)_*(E)$. Since $E$ is the exceptional divisor, we have
$(p_X)_*(E)=0$, an thus $p_X^*(w).\Y_t.E = p_X^*(z).E =0$. Now the
intersection of $\Y_t$ with $E$ is $C$, so we see that $p_X^*(w).C
=0$. To contradict this, we need to prove that $ \Pic(X) \to \Pic(C) $
is injective. This is done as follows.

        Let $G$ be the Grassmannian of all lines in $S$, and $\tilde{X}
\into X \times G$ be the incidence locus.  By choosing a large enough
$n$, we can ensure that the fibres of $p_G:\tilde{X} \to G$ are
irreducible outside a codimension three subset of $G$. This can be
done as follows. Let $R \subset \tilde{X}$ be a subset where the
morphism $p:\tilde{X} \to G$ is not smooth. Then since $p$ is
generically smooth, $R$ is a proper closed subscheme of $\tX$.

\begin{claim}
    $R$ is an irreducible subscheme of $\tilde{X}$.
 \end{claim}
 \begin{proof}
        For simplicity, we write $S=\P^n$. By definition, $(x,L) \in
\tX$ is in $R$ if and only if the the natural map $T_\tX \to p^*T_G$
is not surjective. Note that lines in the dual projective space
$\tilde{\P^n}$ correspond to codimension two linear spaces in $\P^n$.
We will use the same notation to denote a line in $\tilde{\P^n}$ and
the corresponding codimension two linear space in $\P^n$. Thus we see
that $p$ is not smooth at $(x,L)$ if and only if the tangent spaces
$T_{X,x}$ and $T_{L,x}$ do not span the tangent space $T_{\P^n,x}$.
One has the following diagram of vector spaces:
 $$
 \Matrix{
        0 & \to & T_{L,x} & \to & T_{\P^n,x} & \to & N_{L/\P^n} & \to
                & 0 \cr
         &  &  & & \parallel &  & & & \cr
        0 & \to & T_{X,x} & \to & T_{\P^n,x} & \to & N_{X/\P^n} & \to
                & 0 \cr
        }
 $$
         Then $R$ is the locus of pairs $(x,L)$ where $\coker(T_{L,x}
\to N_{X/\P^n,x}) \neq 0$.  So on writing $U=T_{\P^n,x},
E=N_{X/\P^n,x}$, for each $L$, such that $(x,L)\in R$, we have a
codimension $n-2$ subspace $U_L$ of $U$.  We are interested in those
which subspaces $U_L$ whose images under the surjection $U\to E$ of
non-zero co-rank. Write $K=\ker(U\to E)$.  Clearly the composite map
$U_L \to W$ is not surjective if and only if we have $\dim U\cap K
\geq \dim U -\dim W +1$.  Thus we see that the set of such subspaces
$U_L$ is in fact a Schubert subvariety in the Grassmannian of
codimension two subspaces of $U$. It is well know that such varieties
are irreducible.

        Now varying $x\in X$ we see that, for every $x\in X$, we have
a surjection of vector spaces $T_{\P^n,x} \to N_{X/\P^n,x}$ and the
set of pairs $(x,L)\in R$ is a codimension two linear subspace
$T_{L,x} \into T_{\P^n,x}$. Furthermore, $(x,L)\in R$ if and only if
the composite map $T_{L,x} \to N_{X/\P^n,x}$ is not surjective. Thus
$R$ is an irreducible variety.
 \end{proof}

        Now the subscheme $R'$ of $\tX$ where the fibres of the map
$p:\tX \to G$ are of dimension two is in fact a subscheme of $R$. It
is in fact a proper subscheme of $R$. This follows from the fact that
since $n$ is large enough, the general singular fibre of $p$ (which a
complete intersection) is irreducible. This is an easy consequence of
Bertini's Theorem. Then the image of $R'$ in $G$ has  codimension at
least three.

        Now let $Z$ be the open subset of $G$ where the fibres of
$p_G$ are smooth and irreducible. Let $\eta$ be the generic fibre of
$Z$.

 We have a complex
        $$\oplus_D CH^0(p_G^{-1}(D)) \to \Pic(\tilde{X}) \to
                \Pic(\tilde{X}_\eta) $$
 where the sum extends over all codimension one subvarieties $D$ in
$G$. By the choice of $n$ as above, $p_G^{-1}(D) $ is irreducible.
Therefore, the image of $CH^0(p_G^{-1}(D)) \to \Pic(\tilde{X})$ is
contained in the image of $\Pic(G) \to \Pic(\tilde{X})$.  Note that
$\Pic(\tilde{X}) =\Pic(X) \oplus \Z$. Thus $\Pic(X) \to
\Pic(\tilde{X}_\eta)$ is injective.

         Now the following proposition, which seems to be well known
but for which we have been unable to find a convenient reference,
completes the proof.

\begin{propose}\label{specialization}
    Let $\pi: \tilde{X} \to X, f:\tilde{X}\to Z$ be such that  $f$ is
proper with general  fibre irreducible and smooth,  and $X,Z$
irreducible, $X$ smooth projective variety over complex numbers.  Let
$\eta$ be the generic point of $Z$.  Suppose $\Pic(X) \to
\Pic(\tilde{X}_\eta)$ is injective.  Then for $s\in Z$ outside a
countable union of proper closed subvarieties of $Z$, the
``specialization map'' $\Pic(X) \to \Pic(\tilde{X}_s)$ is injective.
 \end{propose}

\begin{proof}
        Recall  that the N\'eron-Severi group, $NS(X)
=\Pic(X)/\Pic^0(X)$, is a finitely generated group. Let $\tau \in
NS(X)$ be a numerical class, then it is well-known (see for instance
\cite{Mumford-Curves})  that there exists smooth projective variety
$\Pic^\tau(X)$, which parametrises line bundles on $X$ of numerical
class $\tau$, and a universal line bundle $P_\tau$ on $X \times
\Pic^\tau(X)$ (the ``Poincare bundle'') with the following property:
for any $\alpha\in \Pic^\tau(X)$ the restriction
 $P_\tau\big|_{X\times \{\alpha\}}$ is the line bundle $\alpha$ on
$X$. The line bundle $P_\tau$ is unique up to a tensoring with a line
bundle pulled back from $\Pic^\tau(X)$. We have then the following
diagram of morphisms:
\def\id{\mathop{\rm id}\nolimits}
        $$\Matrix{
                \tilde{X}\times \Pic^\tau(X) &
                        \mapright{f'=f\times\id} & Z \times
                                        \Pic^\tau(X) \cr
                \mapdown{\pi\times\id} & & \cr
                X\times \Pic^\tau(X) & & \cr
                 }
        $$
 Now recall that $\Pic(X) = \coprod_{\tau\in NS(X)} \Pic^\tau(X)$.
Let $s\in Z$, then the ``specialization map'' $g:\Pic(X) \to
\Pic(\tX_s)$ gives to maps $g_{\tau,s}:\Pic^\tau(X) \to \Pic(\tX_s)$.
Note that $0\in g_{0,s}^{-1}(0)$ for any $s\in Z$. Then to prove
the proposition, it suffices to prove that for $s$ outside a countable
union of proper closed subsets, we have
        $$
        g_{\tau,s}^{-1}(0)  = \cases{\{0\}, &if $\tau = 0$;\cr
                                \emptyset, & otherwise.\cr
                                }
        $$

Let $P'_\tau = (\pi\times\id)^*(P_\tau)$. Then we are interested in
the set of points $Z_{\tau}=\{(s,p)\in Z\times \Pic^\tau(X) \big|
P'_\tau\big|_{\tilde{X}_s\times\{p\}}\} = 0$. We claim that $Z_\tau$ is in
fact a closed subset of $Z\times\Pic^\tau(X)$. This is accomplished
by the following well-known:

\begin{schol}
    Let $X$ be a smooth projective variety, $S$ any irreducible
smooth variety. Let $X \to S$ be smooth proper morphism, and let
${\cal E}$ be a coherent sheaf on $X$ which is flat over $S$. Then
the set of points $s\in S$ where the $H^0(X_s,{\cal E}_s) \neq 0 $ is
a closed subset of $S$.
 \end{schol}
         The assertion of the Scholium is clearly local, so we can
assume that the base is affine. Now the result is an easy consequence
of the existence of a Grothendieck complex for ${\cal E}$ (see
\cite{Mumford}), and will be left to the reader.

        Now we apply the Scholium to the morphism $f':\tX \times
\Pic^\tau(X) \to X\times \Pic^\tau(X)$ and with ${\cal E} = P'_\tau$,
and with ${\cal E}  = P_\tau^{'-1}$. Now we are done because a line
$L$ bundle on an irreducible variety is trivial if and only if $H^0(L)
\neq 0$ and $H^0(L^{-1})\neq 0$. So the set $Z_\tau$ is closed.

  Now the proof breaks up into the above two cases: $\tau=0$ and
$\tau\neq 0$. Let us first dispose of the case $\tau \neq 0$. Since by
the hypothesis of Proposition~\ref{specialization} $\Pic(X) \to
\Pic(\tX_\eta)$ is injective, and since the morphism $f'$ is
generically smooth we see that the subset $Z_\tau$ does not meet
$\{\eta\}\times \Pic^\tau(X)$. Hence, there is a proper closed subset
$D_\tau\subset Z$ such that $Z_\tau \into D_\tau\times \Pic^\tau(X)$.
Hence the image of $Z_\tau$ in $Z$ under the first projection is a
proper closed subset of $Z$.

        The argument for $\tau=0$ is almost identical except that at
the generic point, the line bundle $P'_0$ is trivial, and hence at the
generic point, $Z_0$ is  contained in a subset $\eta\times\Pic^0(X)
\cup D_0\times\Pic^0(X)$, where $D_0$ is a proper closed subset of
$X$.  So to sum up this argument, we have shown that for every class
$\tau$ in $NS(X)$, there is closed subscheme $Z_\tau$ of $Z\times
\Pic^\tau(X)$, such that for $\tau\neq 0$, $Z_\tau\into D_\tau\times
\Pic^\tau(X)$, where $D_\tau\subset Z$ is a proper closed subscheme.
For $\tau=0$, the subscheme $Z_0$ (at the generic point) is contained
in a subscheme of the form $\eta\times\Pic^0(X)\cup D_0\times
\Pic^0(X)$, and again $D_0$ is a proper subscheme of $Z$.  Since
$NS(X)$ is finitely generated, it is countable, so by removing a
countable union of closed subsets  $D_\tau$ where $\tau\in NS(X)$, we
can ensure that the specialization map is injective.
 \end{proof}

\begin{rmk}
            In a personal communication, Nori has pointed out to us
that  his conjecture is not valid for Chow groups over $\Z$, in
other words one has to work with Chow groups tensor $\Q$. This can
been seen as follows. Suppose $X$ is a smooth projective threefold,
$L$ an ample line bundle on $X$, assume further that $H^1(X,\O_X) =
0$. Let $n$ be any positive integer. Then the zero scheme $Y$ of a
general section of $H^0(X,L^n)$ has the property that
$H^1(\Y_{K},\O_{\Y_K}) =0 $. Then by a well-known result of Roitman
(see \cite{Roitman}), we know that $CH^2(\Y_{\bar K})$ is torsion
free, because the Albanese variety of $\Y_{\bar K}=0$. So that the
restriction $CH^2(X_{K}) \to CH^2(\Y_{\bar K})$ obliterates all the
torsion in $CH^2(X_{K})$, consequently if $CH^2(X_{K})$ has
non-trivial torsion then the restriction map cannot be injective.

        It is possible to write down examples of threefolds which have
the above properties. A well studied example being that of a smooth
general quartic hypersurface in $\P^4$.  It is shown in \cite{Bloch}
that group of codimension two cycles algebraically equivalent to zero
modulo rational equivalence, denoted by  $A^2(X)$, is isomorphic as a
group to the intermediate Jacobian $J^2(X)$. The latter group is a
complex torus (in fact an Abelian variety) and hence has non-trivial
torsion. Moreover, it is standard that $A^2(X) \subset CH^2(X)$ (see
\cite{Fulton}). A similar example can be found in \cite{Murre}.
 \end{rmk}

 \end{proof}

 \bibliographystyle{plain}

\end{document}